\input amstex
\documentstyle{amsppt} \nologo
\loadbold

\let\bs\boldsymbol
\hsize=5.75truein
\vsize=8.75truein 
\hcorrection{.25truein}
\loadeusm \let\scr\eusm
\loadeurm 
 at 17truept
\font\Rrm=cmr17 
\font\Iit=cmmi12 scaled \magstep2
\define\pre#1#2{{_{#1}{#2}}} 
\define\cind{\text{\it c-\/\rm Ind}} 
\define\Aut#1#2{\text{\rm Aut}_{#1}(#2)}
\define\End#1#2{\text{\rm End}_{#1}(#2)}
\define\GL#1#2{\roman{GL}_{#1}(#2)}
\define\M#1#2{\roman M_{#1}(#2)}
\define\bsi{\bs\iota} 
\define\upr#1#2{{}^{#1\!}{#2}} 
\define\Ind{\text{\rm Ind}}
\let\ge\geqslant
\let\le\leqslant
\let\vP\varPi
\let\vt\vartheta
\let\vT\varTheta 
\topmatter 
\title\nofrills \Rrm 
Intertwining of simple characters in GL({\Iit n\/})
\endtitle 
\rightheadtext{Simple characters} 
\author 
Colin J. Bushnell and Guy Henniart 
\endauthor 
\leftheadtext{C.J. Bushnell and G. Henniart}
\affil 
King's College London and Universit\'e de Paris-Sud 
\endaffil 
\address 
King's College London, Department of Mathematics, Strand, London WC2R 2LS, UK. 
\endaddress
\email 
colin.bushnell\@kcl.ac.uk 
\endemail
\address 
Institut Universitaire de France et Universit\'e de Paris-Sud, Laboratoire de Math\'ematiques d'Orsay,
Orsay Cedex, F-91405; CNRS, Orsay cedex, F-91405. 
\endaddress 
\email 
Guy.Henniart\@math.u-psud.fr 
\endemail 
\date July 2011 \enddate 
\thanks 
This paper was conceived while the first-named author was visiting, and partially supported by, Universit\'e de Paris-Sud. 
\endthanks 
\abstract 
Let $F$ be a non-Archimedean local field and let $G$ be the general linear group $G = \text{\rm GL}_n(F)$. Let $\theta_1$, $\theta_2$ be simple characters in $G$. We show that $\theta_1$ intertwines with $\theta_2$ if and only if $\theta_1$ is endo-equivalent to $\theta_2$. We also show that any simple character in $G$ is a $G$-type.  \endabstract 
\endtopmatter 
\document \nopagenumbers \baselineskip=13pt \parskip=3pt plus 1pt
Let $F$ be a non-Archimedean local field and let $G = \GL nF$, for some $n\ge1$. Following \cite{1}, the category $\roman{Rep}\,G$ of smooth complex representations of $G$ decomposes as a direct sum of indecomposable blocks, 
$$
\roman{Rep}\,G = \coprod_{\frak s\in \scr B(G)}\,\roman{Rep}_\frak s\,G, 
$$ 
indexed by a certain set $\scr B(G)$. Let $\frak S$ be a finite subset of $\scr B(G)$. As in \cite{7}, an $\frak S$-type in $G$ is an irreducible smooth representation $\rho$, of some compact open subgroup of $G$, with the property that an irreducible smooth representation of $G$ contains $\rho$ if and only if it lies in $\roman{Rep}_\frak s\,G$, for some $\frak s\in \frak S$. 
\par 
One knows \cite{8} how to construct an $\{\frak s\}$-type in $G$, for any $\frak s\in \scr B(G)$. Those types are all built from {\it simple characters\/} in groups $\GL mF$, in the sense of \cite{6}, for various integers $m\le n$. Here, we return to the simple characters themselves and prove: 
\proclaim{Type Theorem} 
Let $G = \GL nF$, for some $n\ge 1$, and let $\theta$ be a simple character in $G$. The character $\theta$ is then an $\frak S_\theta$-type in $G$, for some finite subset $\frak S_\theta$ of $\scr B(G)$. 
\endproclaim 
The proof, and a description of $\frak S_\theta$, are given in \S4 below. 
\par 
We use the Type Theorem to prove a powerful result concerning the intertwining properties of simple characters. As part of the definition, a simple character in $G = \GL nF$ is attached, in an invariant manner, to a hereditary order in the matrix algebra $A = \M nF$. A cornerstone of the theory is the fact ((3.5.11) of \cite{6}) that two simple characters in $G$, attached to the same order and which intertwine in $G$, are actually conjugate. This result is taken one further level in \cite{4}. There is a canonical procedure for transferring simple characters between hereditary orders, in possibly different matrix algebras. Given two simple characters $\theta_i$ in $\GL{n_i}F$, attached to hereditary orders $\frak a_i$, one can find an integer $n$ and a hereditary order $\frak a$ in $\M nF$ to which both characters may be transferred. If the transferred characters are conjugate in $\GL nF$, one says they are {\it endo-equivalent.} One knows that endo-equivalence is an equivalence relation on the class of simple characters in all general linear groups. Moreover, endo-equivalent characters attached to the same order are necessarily conjugate. Here, we consider a general pair of simple characters in one group. We prove: 
\proclaim{Intertwining Theorem} 
Let $\theta_1$, $\theta_2$ be simple characters in $G = \GL nF$. The characters $\theta_1$, $\theta_2$ intertwine in $G$ if and only if they are endo-equivalent. 
\endproclaim 
To the specialist in the area, these results provide clear and satisfying conclusions to several lines of development, but the non-specialist may wish for more motivation. This is first provided by our examination \cite{5} of the congruence properties of the local Langlands correspondence, where the Intertwining Theorem provides a crucial step in the argument: see \cite{5} 4.3 Lemma. 
\par
The results also give a framework in which to investigate representations in more general settings. There is a fully functional theory of simple characters and endo-equivalence spanning all inner forms $\GL mD$ of $\GL nF$, where $D$ is a finite-dimensional central $F$-division algebra \cite{2}, \cite{9}, \cite{10}. However, certain new structures come into play, and one is led to ask how these are reflected or clarified in analogues of our results. In a different direction, one may consider smooth representations of $\GL nF$ over fields of positive characteristic $\ell$. Provided $\ell$ is not the residual characteristic of $F$, one has an identical general theory of simple characters. One deduces readily that the Intertwining Theorem holds unchanged. However, as Vincent S\'echerre reminds us, a simple character need not be a type in this situation. 
\head 
1. A review of simple characters 
\endhead 
Let $\frak o_F$ be the discrete valuation ring in $F$ and $\frak p_F$ the maximal ideal of $\frak o_F$. We choose a smooth character $\psi$ of $F$ which is of level one, in the sense that $\roman{Ker}\,\psi$ contains $\frak p_F$ but not $\frak o_F$. 
\par 
Let $V$ be an $F$-vector space of finite dimension and set $A = \End FV$, $G = \Aut FV$. Let $\frak a$ be a hereditary $\frak o_F$-order in $A$, with Jacobson radical $\frak p_\frak a$. A {\it simple character in $G$,} attached to $\frak a$, is one of the following objects. The {\it trivial\/} simple character attached to $\frak a$ is the trivial character of the group $U^1_\frak a = 1{+}\frak p_\frak a$: we denote this $1^1_\frak a$. 
\par 
To define a {\it non-trivial\/} simple character attached to $\frak a$, we recall briefly the definition \cite{6} (1.5.5) of a simple stratum $[\frak a,l,0,\beta]$ in $A$. First, $\beta\in G$ and the algebra $E = F[\beta]$ is a field. The hereditary order $\frak a$ is {\it $E$-pure,} in that $x^{-1}\frak a x =\frak a$ for $x\in E^\times$. The integer $l$ is positive and given by $\beta^{-1}\frak a = \frak p_\frak a^l$. The quadruple $[\frak a,l,0,\beta]$ is then a simple stratum in $A$ if $\beta$ satisfies a technical condition ``$k_0(\beta,\frak a)<0$'' {\it loc\. cit.} Since we will not use this directly, we say no more of it. 
\par
Following the recipes of \cite{6} 3.1, the simple stratum $[\frak a,\beta] = [\frak a,l,0,\beta]$ defines open subgroups $H^1(\beta,\frak a)\subset J^1(\beta,\frak a) \subset J^0(\beta,\frak a)$ of the unit group $U_\frak a = \frak a^\times$, such that $J^1(\beta,\frak a) = J^0(\beta,\frak a)\cap U^1_\frak a$. The choice of $\psi$ then gives rise to a finite set $\scr C(\frak a,\beta,\psi)$ of smooth characters of $H^1(\beta,\frak a)$, called {\it simple characters:\/} see \cite{6} 3.2 for the full definition. The choice of $\psi$ is essentially irrelevant, so we treat it as fixed and henceforth omit it from the notation. 
\par 
We recall a fundamental property of simple characters attached to a fixed hereditary order \cite{6} (3.5.11). 
\proclaim{Intertwining implies conjugacy} 
For $i=1,2$, let $[\frak a,\beta_i]$ be a simple stratum in $A$ and let $\theta_i\in \scr C(\frak a,\beta_i)$. If the characters $\theta_1$, $\theta_2$ intertwine in $G$ then they are conjugate by an element of $U_\frak a$. 
\endproclaim 
We shall also need systems of {\it transfer maps.} Let $[\frak a,l,0,\beta]$ be a simple stratum in $A$, as before. Suppose we have another $F$-vector space $V'$ of finite dimension, an $F$-embedding $\bsi:F[\beta] \to A' = \End F{V'}$, and an $F[\bsi \beta]$-pure hereditary order $\frak a'$ in $A'$: any two such embeddings $\bsi$ are $U_{\frak a'}$-conjugate, so we are justified in omitting $\bsi$ from the notation. There is a unique integer $l'$ such that $[\frak a',l',0,\beta]$ is a simple stratum in $A'$. There is a canonical bijection 
$$
\tau^\beta_{\frak a,\frak a'}: \scr C(\frak a,\beta) @>{\ \ \approx\ \ }>> \scr C(\frak a',\beta). 
\tag 1.1 
$$ 
This family of maps is transitive with respect to the orders: in the obvious notation, we have 
$$
\tau^\beta_{\frak a,\frak a''} = \tau^\beta_{\frak a',\frak a''}\circ \tau^\beta_{\frak a,\frak a'}. 
$$ 
Full details may be found in \cite{6} section 3.6 and \cite{4} section 8. 
\proclaim{Lemma 1} 
For $j=1,2$, let $[\frak a_j,l_j,0,\beta_j]$ be a simple stratum in $A_j = \End F{V_j}$. 
\roster 
\item 
There exists a finite-dimensional $F$-vector space $V$, a hereditary order $\frak a$ in $A = \End FV$ and a pair of $F$-embeddings $\bsi_j:F[\beta_j]\to A$, such that $\frak a$ is $F[\bsi_j\beta_j]$-pure, for $j=1,2$. 
\item 
Let $\theta_j\in \scr C(\frak a_j,\beta_j)$. The following are equivalent: 
\itemitem{\rm (a)} 
There exists a system $(V,\frak a,\bsi_j)$, as in \rom{(1),} such that $\tau^{\beta_1}_{\frak a_1,\frak a}\theta_1$ intertwines with $\tau^{\beta_2}_{\frak a_2,\frak a}\theta_2$ in $G = \Aut FV$. 
\itemitem{\rm (b)} 
For any system $(V,\frak a,\bsi_j)$, as in \rom{(1),} the character $\tau^{\beta_1}_{\frak a_1,\frak a}\theta_1$ intertwines with $\tau^{\beta_2}_{\frak a_2,\frak a}\theta_2$ in $G = \Aut FV$. 
\endroster 
\endproclaim 
\demo{Proof} 
Part (1) is elementary. If (2)(a) holds, then $[F[\beta_1]{:}F] = [F[\beta_2]{:}F]$ by \cite{6} (3.5.1) and the intertwining implies conjugacy property. Part (2) is then given by Theorem 8.7 of \cite{4}. \qed 
\enddemo 
Developing this theme, if the (non-trivial) simple characters $\theta_j$ of Lemma 1(2) satisfy condition (a), we say they are {\it endo-equivalent.} In particular, in the context of (1.1), $\theta$ is endo-equivalent to $\tau^\beta_{\frak a,\frak a'}\theta$. Further, two endo-equivalent simple characters attached to the same order intertwine, and so are conjugate. It follows that endo-equivalence is indeed an equivalence relation on the class of non-trivial simple characters {\it cf\.} \cite{4} 8.10. 
\par 
It is convenient to extend this framework to include the trivial simple characters. We set $\tau_{\frak a,\frak a'}1^1_\frak a = 1^1_{\frak a'}$ and deem that any two trivial simple characters are endo-equivalent. Any two such characters in the same group intertwine, so the approach is consistent with the main case. Moreover, a trivial simple character can never intertwine with a non-trivial one: this follows from \cite{3} Theorem 1 and \cite{6} (2.6.2). 
\head 
2. Heisenberg extensions  
\endhead 
Let $\theta$ be a simple character in $G = \Aut FV$, attached to a hereditary order $\frak a$ in $A = \End FV$. Thus $\theta$ is a character of an open subgroup $H^1_\theta$ of $U^1_\frak a$. Let $J^0_\theta$ denote the $U_\frak a$-normalizer of $\theta$ and put $J^1_\theta = J^0_\theta\cap U^1_\frak a$. If $\theta$ is non-trivial, we choose  a simple stratum $[\frak a,\beta]$ in $A$ such that $\theta\in \scr C(\frak a,\beta)$. We then get $J^k_\theta = J^k(\beta,\frak a)$ and $H^1_\theta = H^1(\beta,\frak a)$, in the notation of \S1. If $\theta$ is the trivial simple character $1^1_\frak a$ attached to $\frak a$, we have $H^1_\theta = J^1_\theta = U^1_\frak a$ and $J^0_\theta = U_\frak a$. 
\par 
With $E = F[\beta]$ (if $\theta$ is non-trivial) or $F$ (otherwise), let $B  = \End EV$ be the $A$-centralizer of $E$ and set $\frak b = \frak a\cap B$. Thus $\frak b$ is a hereditary $\frak o_E$-order in $B$ with radical $\frak q = \frak p_\frak a\cap B$. We then have $J^0_\theta = J^1_\theta U_\frak b$ and $J^1_\theta\cap U_\frak b = U^1_\frak b$. 
\par
Let $\eta = \eta_\theta$ be the unique irreducible representation of $J^1_\theta$ which contains $\theta$ \cite{6} (5.1.1). Thus $\eta|_{H^1_\theta}$ is a multiple of $\theta$. Let $\scr R^0(\theta)$ be the set of equivalence classes of irreducible representations of $J^0_\theta$ which contain $\theta$.  Let $\scr H^0(\theta)$ be the set of $\kappa\in \scr R^0(\theta)$ with the following two properties. First, $\kappa|_{J^1_\theta} \cong \eta_\theta$. Second, $\kappa$ is intertwined by every element of $G$ which intertwines $\theta$. In the language of \cite{6}, $\scr H^0(\theta)$ consists of the ``$\beta$-extensions'' of $\eta_\theta$ and is non-empty \cite{6} (5.2.2). 
\par 
In particular, $\scr R^0(1^1_\frak b)$ is the set of equivalence classes of irreducible representations of $U_\frak b$ trivial on $U^1_\frak b$. For $\sigma\in \scr R^0(1^1_\frak b)$, there is a unique irreducible representation $\sigma_\theta$ of $J^0_\theta$ which agrees with $\sigma$ on $U_\frak b$ and is trivial on $J^1_\theta$. 
\proclaim{Lemma 2} 
Let $\kappa\in \scr H^0(\theta)$, $\sigma\in \scr R^0(1^1_\frak b)$. The representation $\kappa\otimes \sigma_\theta$ of $J^0_\theta$ is irreducible, and lies in $\scr R^0(\theta)$. For any $\kappa\in \scr H^0(\theta)$, the map 
$$
\align 
\scr R^0(1^1_\frak b) &\longrightarrow \scr R^0(\theta), \\ \sigma&\longmapsto \kappa\otimes \sigma_\theta, 
\endalign 
$$ 
is a bijection. 
\endproclaim 
\demo{Proof} 
The restriction of $\kappa\otimes\sigma_\theta$ to $H^1_\theta$ is surely a multiple of $\theta$. The other assertions are given by \cite{5} 1.5 Proposition. \qed 
\enddemo 
\head 
3. Residually cuspidal representations 
\endhead 
We continue in the same situation. Let $\Bbbk_E$ denote the residue field of $E$. The group $J_\theta^0/J_\theta^1$ takes the form 
$$
J_\theta^0/J_\theta^1 \cong U_\frak b/U^1_\frak b \cong \prod_{i=1}^r \GL{m_i}{\Bbbk_E}, 
$$
for integers $r, m_i\ge1$ such that $\sum_{1\le i\le r} m_i = n/[E{:}F]$. In particular, $U_\frak b/U^1_\frak b$ is the group of rational points of a connected reductive $\Bbbk_E$-group. We fix $\kappa\in \scr H^0(\theta)$. If $\lambda\in \scr R^0(\theta)$ then, by Lemma 2, $\lambda\cong \kappa\otimes \sigma_\theta$ where $\sigma$ is the inflation of a uniquely determined irreducible representation $\tilde\sigma$ of $U_\frak b/U^1_\frak b$. We say that $\lambda$ is {\it residually cuspidal\/} if the representation $\tilde\sigma$ is cuspidal. The representation $\kappa$ is uniquely determined, up to tensoring with a character of the form $(\phi\circ\det_B)_\theta$, where $\phi$ is a character of $U_E$ trivial on $U^1_E$ \cite{6} (5.2.2), so this property of $\lambda$ does not depend on the choice of $\kappa$. We denote by $\scr R^0_c(\theta)$ the subset of residually cuspidal elements of $\scr R^0(\theta)$. 
\proclaim{Proposition 1} 
Let $\theta$ be a simple character in $G$, and let $[\frak a,\beta]$ be a simple stratum in $A$ such that $\theta\in \scr C(\frak a,\beta)$. Let $E$ denote the field $F[\beta]$. 
\roster 
\item 
Let $\frak a'$ be an $E$-pure hereditary $\frak o_F$-order in $A$, containing $\frak a$. Let $\theta' = \tau^\beta_{\frak a,\frak a'}\theta$, and let $\lambda\in \scr R^0_c(\theta)$. An irreducible representation $\pi$ of\/ $G$ containing $\lambda$ then contains some element of\/ $\scr R^0(\theta')$. 
\item 
Suppose that $\lambda\in \scr R^0(\theta)$ is not residually cuspidal. There exists an $E$-pure hereditary $\frak o_F$-order $\frak a''$ in $A$, with $\frak a''\varsubsetneq \frak a$, and an element $\lambda'' \in \scr R^0_c(\theta'')$, where $\theta'' = \tau^\beta_{\frak a,\frak a''}\theta$, with the following property: any irreducible representation of $G$ containing $\lambda$ also contains $\lambda''$. 
\endroster 
\endproclaim 
\demo{Proof} 
All assertions follow from (8.3.5) Proposition of \cite{6}. \qed 
\enddemo 
\remark{Remark 1} 
Proposition 1 applies equally when $\theta$ is a trivial simple character, as noted in \cite{6}, following (8.3.5). 
\endremark 
The simple characters $\theta'$, $\theta''$ of Proposition 1 are both endo-equivalent to $\theta$. 
\par 
We say that a simple stratum $[\frak a,\beta]$ in $A$ is {\it m-simple\/} if $\frak a$ is maximal among $F[\beta]$-pure hereditary $\frak o_F$-orders in $A$. We say that a simple character $\theta$ is {\it m-simple\/} if $\theta\in \scr C(\frak a,\beta)$, where $[\frak a,\beta]$ is m-simple. (This depends on $\theta$, not the choice of $[\frak a,\beta]$.) Similarly for trivial characters. 
\proclaim{Proposition 2} 
Let $\lambda\in \scr R^0(\theta)$. The following are equivalent: 
\roster 
\item $\theta$ is m-simple and $\lambda$ is residually cuspidal; 
\item $\lambda$ is contained in some irreducible cuspidal representation of $G$; 
\item any irreducible representation of $G$ containing $\lambda$ is cuspidal. 
\endroster 
\endproclaim 
\demo{Proof} 
The equivalence of (2) and (3) is \cite{6} (6.2.1, 6.2.2). The implication $(1)\Rightarrow (2)$  is \cite {5} (6.2.3). For the converse, suppose that either $\theta$ is not m-simple or that $\lambda$ is not residually cuspidal. Let $\pi$ be an irreducible representation of $G$ containing $\lambda$. In either case, part (2) of Proposition 1 implies the existence of the following objects: 
\roster 
\item a non-maximal $E$-pure hereditary order $\frak a'$ in $A$, 
\item a simple character $\theta'$ attached to $\frak a'$ and endo-equivalent to $\theta$, 
\item a representation $\lambda'\in \scr R^0_c(\theta')$ occurring in $\pi$. 
\endroster 
The representation $\pi$ is then not cuspidal, by \cite{6} (8.3.3 or 7.3.16). \qed 
\enddemo 
\proclaim{Corollary 1} 
An irreducible cuspidal representation $\pi$ of\/ $G$ contains exactly one conjugacy class of simple characters $\theta$, and all of those characters are m-simple. 
\endproclaim 
\demo{Proof} 
This follows from Proposition 2 and \cite{6} (6.2.4). \qed 
\enddemo 
So, if $\pi$ is an irreducible cuspidal representation of $G$, all simple characters contained in $\pi$ belong to the same endo-equivalence class, which we denote $\vt(\pi)$. 
\head 
4. The Type Theorem 
\endhead 
Let $\frak a$ be a hereditary $\frak o_F$-order in $A = \End FV$, with Jacobson radical $\frak p_\frak a$. Thus 
$$
\frak a/\frak p_\frak a \cong \prod_{i=1}^r \M{n_i}{\Bbbk_F}, 
$$
for positive integers $n_i$ with sum $n$. Let $M_\frak a$ be an $F$-Levi subgroup of $G$ such that $M_\frak a\cong \prod_{1\le i\le r} \GL{n_i}F$. The group $M_\frak a$ is determined uniquely, up to conjugation in $G$. If $M$ is an $F$-Levi subgroup of $G$, we say that $M$ is {\it subordinate to\/} $\frak a$ if $M$ is $G$-conjugate to a Levi subgroup of $M_\frak a$. 
\par
We recall some further definitions. A {\it cuspidal datum in $G$} is a pair $(M,\sigma)$, where $M$ is a Levi subgroup of $G$ and $\sigma$ is an irreducible cuspidal representation of $M$. The set of such data carries the equivalence relation ``$G$-inertial equivalence'', as in \cite{7} \S1. The set of equivalence classes for this relation will be denoted $\scr B(G)$. 
\par 
If $\pi$ is an irreducible smooth representation of $G$, there is a cuspidal datum $(M,\sigma)$ in $G$ and a parabolic subgroup $P$ of $G$, with Levi component $M$, such that $\pi$ is equivalent to a subquotient of the induced representation $\Ind_P^G\,\sigma$. The inertial equivalence class of $(M,\sigma)$ is thereby uniquely determined: we call it the {\it inertial support\/} of $\pi$ and denote it $\scr I(\pi)$. If $\frak S$ is a finite subset of $\scr B(G)$, an $\frak S$-type in $G$ is a pair $(K,\rho)$, where $K$ is a compact open subgroup of $G$ and $\rho$ is an irreducible smooth representation of $K$ such that, if $\pi$ is an irreducible smooth representation of $G$, then $\pi$ contains $\rho$ if and only if $\scr I(\pi)\in \frak S$ \cite{7} 4.1, 4.2. 
\par 
Let $\theta$ be a (possibly trivial) simple character in $G$, attached to the hereditary order $\frak a$. Let $\vT$ denote the endo-equivalence class of $\theta$. 
\proclaim{Definition} 
Let $\frak s\in \scr B(G)$ be the $G$-inertial equivalence class of $(M,\sigma)$, where 
$$
M \cong \prod_{j=1}^s \GL{m_j}F, \qquad \sigma = \bigotimes_{j=1}^s \sigma_j, 
$$
and $\sigma_j$ is an irreducible cuspidal representation of $\GL{m_j}F$. We say that $\frak s$ is\/ {\rm subordinate to} $\theta$ if $M$ is subordinate to $\frak a$ and $\vt(\sigma_j) = \vT$, for all $j$. 
\endproclaim 
We prove the following version of the Type Theorem. 
\proclaim{Theorem 3} 
Let $V$ be a finite-dimensional $F$-vector space, and let $G$ denote the group $\Aut FV$. Let $\theta$ be a simple character in $G$. Let $\frak S_\theta$ be the set of $\frak s\in \scr B(G)$ that are subordinate to $\theta$. The character $\theta$ is then an $\frak S_\theta$-type in $G$. 
\endproclaim 
\demo{Proof} 
We have to show that an irreducible representation $\pi$ of $G$ contains $\theta$ if and only if the inertial support of $\pi$ is an element of $\frak S_\theta$. We assume that $\theta$ is non-trivial: the proof for trivial simple characters is parallel but easier, so we omit it. We choose a simple stratum $[\frak a,\beta]$ such that $\theta\in \scr C(\frak a,\beta)$ and set $E = F[\beta]$. 
\par 
We start in a slightly more general situation, with a cuspidal datum $\frak s$ of the form $(M,\sigma)$ such that 
$$ 
M\cong \prod_{k=1}^s\GL{n_k}F, \quad \sigma = \bigotimes_{k=1}^s \sigma_k, 
\tag 4.1 
$$ 
for various integers $n_k\ge1$, and $\vt(\sigma_k) = \vT$ for all $k$. Replacing $M$ by a $G$-conjugate and each $\sigma_k$ by an equivalent representation, we can assume we are in the following situation. First, $M$ is the $G$-stabilizer of a decomposition $V = \bigoplus_{1\le k\le s} V_k$, in which the $V_k$ are non-zero $E$-subspaces of $V$. Second, each $\sigma_k$ contains a simple character $\theta_k\in \scr C(\frak a_k,\beta)$, endo-equivalent to $\theta$, for some simple stratum $[\frak a_k,\beta]$ in $\End F{V_k}$. Observe that, by Corollary 1, each $\theta_k$ is m-simple, so $J^0_{\theta_k}/J^1_{\theta_k} \cong \M {n_k/[E:F]}{\Bbbk_E}$. 
\par 
We may impose a further normalization. We suppose given an $E$-pure hereditary $\frak o_F$-order $\frak A$ in $A = \End FV$ such that $\frak A/\frak p_\frak A \cong \prod_k\GL{n_k}{\Bbbk_F}$, the integers $n_k$ being as in (4.1). There is then an integer $N>0$ such that $[\frak A,N,0,\beta]$ is a simple stratum in $A$. Let $\theta_\frak A = \tau^\beta_{\frak a,\frak A}\theta \in \scr C(\frak A,\beta)$. In particular,  $\theta_\frak A$ is endo-equivalent to $\theta$. Theorem 7.2 of \cite{8} gives an $\frak s$-type in $G$ of the form $(J^0(\beta,\frak A),\lambda_\frak s)$, where $\lambda_\frak s\in \scr R^0_c(\theta_\frak A)$. 
\remark{Remark 2} 
To be more precise, the construction in \cite{8} 7.2 yields an $\frak s$-type $(K,\tau)$ where, in our notation, $H^1(\beta,\frak A) \subset K \subset J^0(\beta,\frak A)$. The representation of $J^0(\beta,\frak A)$ induced by $\tau$ is our $\lambda_\frak s$. 
\endremark 
Let $\frak s \in \frak S_\theta$. Thus $\frak s$ is subordinate to $\theta$ and we may therefore choose $\frak A\subset \frak a$. Let $\pi$ be an irreducible representation of $G$ of inertial support $\frak s$. By definition, $\pi$ contains $\lambda_\frak s$. By Proposition 1, $\pi$ contains a simple character $\theta'\in \scr C(\frak a,\beta)$ which is endo-equivalent to $\theta_\frak A$. It follows that $\theta'$ is endo-equivalent to $\theta$, and hence $G$-conjugate to $\theta$. In particular, $\pi$ contains $\theta$. 
\par
Conversely, let $\pi$ be an irreducible representation of $G$ which contains $\theta$. Proposition 1(2) gives an $E$-pure hereditary order $\frak a'\subset \frak a$, a simple character $\theta'\in \scr C(\frak a',\beta)$ and a representation $\lambda\in \scr R^0_c(\theta')$ which occurs in $\pi$. Comparing with Theorem 7.2 of \cite{8} again, we see that $\lambda$ is an $\frak s$-type in $G$, for some $\frak s\in \frak S_\theta$. Consequently, the inertial support of $\pi$ is an element of $\frak S_\theta$, as required. \qed 
\enddemo 
\head 
5. The Intertwining Theorem 
\endhead 
We prove the Intertwining Theorem. Let $G = \GL nF$, $A = \M nF$. Let $\theta_1$, $\theta_2$ be simple characters in $G$, with endo-classes $\vT_1$, $\vT_2$ respectively. 
\par
Suppose first that $\vT_1 = \vT_2$. We have to show that $\theta_1$ intertwines with $\theta_2$ in $G$. If the $\theta_i$ are trivial, this is clear, so suppose otherwise. Choose simple strata $[\frak a_i,\beta_i]$ in $A$ such that $\theta_i \in \scr C(\frak a_i,\beta_i)$ and put $E_i = F[\beta_i]$. The relation $\vT_1 = \vT_2$ implies that the field extensions $E_i/F$ have the same ramification indices and the same residue class degrees \cite{4} (8.11). So, there exists an $E_2$-pure hereditary order $\frak a$ in $A$  which is isomorphic to $\frak a_1$. Set $\theta = \tau^{\beta_2}_{\frak a_2,\frak a}\theta_2$. According to Lemma 1, the simple characters $\theta_1$, $\theta$ must intertwine in $G$ and hence be $G$-conjugate, say $\theta = \theta_1^g$, for some $g\in G$. By \cite{6} (3.6.1), the characters $\theta_2$, $\theta$ agree on $H^1(\beta_2,\frak a_2) \cap H^1(\beta_2,\frak a)$. Therefore $g$ intertwines $\theta_1$ with $\theta_2$, as required. 
\par 
For the converse, suppose that $\theta_1$ intertwines with $\theta_2$. Abbreviating $H_i = H^1(\beta_i,\frak a_i)$, this hypothesis implies the existence of a non-trivial $G$ homomorphism 
$$
\cind_{H_1}^G\,\theta_1 \longrightarrow \cind_{H_2}^G\,\theta_2. 
\tag 5.1 
$$ 
Frobenius Reciprocity, for compact induction from open subgroups, implies that the space $\vP_i = \cind_{H_i}^G\,\theta_i$ is generated over $G$ by its $\theta_i$-vectors. Since $\theta_i$ is a type in $G$ (Theorem 3), every irreducible $G$-subquotient of $\vP_i$ contains $\theta_i$ \cite{7} 4.1. The existence of the non-trivial map (5.1) implies there exists an irreducible representation $\pi$ of $G$ containing both $\theta_i$. If the inertial support of $\pi$ is of the form $\big(\prod_k\GL{n_k}F, \bigotimes_k \sigma_k\big)$ then, by Theorem 3 again,  $\vt(\sigma_k) = \vT_1 = \vT_2$, for all $k$. \qed 
\smallskip 
Since endo-equivalence is an equivalence relation, the Intertwining Theorem implies that simple characters, in a fixed group, exhibit the following surprising property. 
\proclaim{Corollary 2} 
Let $\theta_1$, $\theta_2$, $\theta_3$ be simple characters in $G = \GL nF$. If $\theta_1$ intertwines with $\theta_2$ and $\theta_2$ intertwines with $\theta_3$, then $\theta_1$ intertwines with $\theta_3$. 
\endproclaim 
 
\enddocument \bye